\documentclass{amsart}
\usepackage{amssymb}
\usepackage[cp1250]{inputenc}

\newtheorem{theorem}{Theorem}[section]
\newtheorem{lemma}[theorem]{Lemma}
\newtheorem{proposition}[theorem]{Proposition}
\newtheorem{example}[theorem]{Example}

\numberwithin{equation}{section}

\begin{document}

\title{Recurrence relations for characters of affine Lie algebra $A_{\ell}^{(1)}$ }

\author{Miroslav Jerković}
\address{Faculty of Chemical Engineering and Technology, University of Zagreb, Marulićev trg 19, Zagreb, Croatia}
\curraddr{}
\email{mjerkov@fkit.hr}
\thanks{}
\subjclass[2000]{Primary 17B67; Secondary 17B69, 05A19.\\ \indent Partially supported by the Ministry of Science and Technology of the Republic of Croatia, Project ID 037-0372794-2806}
\keywords{}
\date{}
\dedicatory{}

\begin{abstract} By using the known description of combinatorial bases for Feigin-Stoyanovsky's type subspaces of standard modules for affine Lie algebra $\mathfrak{sl}(l+1,\mathbb{C})^{\widetilde{}}$, as well as certain intertwining operators between standard modules, we obtain exact sequences of Feigin-Stoyanovsky's type subspaces at fixed level $k$. This directly leads to systems of recurrence relations for formal characters of those subspaces. \end{abstract}

\maketitle

\section{Introduction}

In a series of papers Lepowsky and Wilson have obtained a Lie theoretic proof of famous Rogers-Ramanujan partition identities (cf. \cite{lw1,lw2}). They used the fact that the product side of these identities arise in principally specialized characters for all level $3$ standard representations of affine Lie algebra $\mathfrak{sl}(2, \mathbb{C})^{\widetilde{}}$, and to obtain the sum side of identities they constructed combinatorial bases for these $\mathfrak{sl}(2, \mathbb{C})^{\widetilde{}}$-modules. The key ingredient in their construction is a twisted vertex operator construction of fundamental $\mathfrak{sl}(2, \mathbb{C})^{\widetilde{}}$-modules and the corresponding vertex operator relations for higher level modules.

By using Lepowsky-Wilson's approach with untwisted vertex operators, Lepowsky and Primc obtained in \cite{lp} new character formulas for all standard $\mathfrak{sl}(2, \mathbb{C})^{\widetilde{}}$-modules. Feigin and Stoyanovsky in \cite{FS} gave another proof and combinatorial interpretation of these character formulas, one of several new ingredients in the proof are the so-called principal subspaces, later also named Feigin-Stoyanovsky's subspaces. By using intertwining operators Georgiev (cf. \cite{G}) extended Feigin-Stoyanovsky's approach and obtained combinatorial bases of principal subspaces of standard $\mathfrak{sl}(\ell + 1, \mathbb{C})^{\widetilde{}}$-modules together with the corresponding character formulas. In a similar way Capparelli, Lepowsky and Milas in \cite{CLM1,CLM2} use the theory of intertwining operators for vertex operator algebras in order to obtain the recurrence relations for characters of principal subspaces of all arbitrary fixed level standard $\mathfrak{sl}(2, \mathbb{C})^{\widetilde{}}$-modules. It turned out that this recurrence relations are precisely the known recursion formulas of Rogers and Selberg, already solved by G. Andrews while working on Gordon identities (cf. \cite{A1,A2}). The Capparelli-Lepowsky-Milas approach was further investigated by Calinescu in \cite{Ca1} in order to construct the exact sequences of principal subspaces of basic $\mathfrak{sl}(l+1,\mathbb{C})^{\widetilde{}}$-modules and thus acquire a recurrence system for characters of these subspaces. Furthermore, in \cite{Ca2} Calinescu used the exact sequence method to obtain systems of recurrences for characters of all principal subspaces of arbitrary level standard $\mathfrak{sl}(3,\mathbb{C})^{\widetilde{}}$-modules and the corresponding characters for some classes of principal subspaces. Finally, this approach was further developed in \cite{CaLM1,CaLM2} where presentations of principal subspaces of all standard $\mathfrak{sl}(l+1,\mathbb{C})^{\widetilde{}}$-modules are given.

Another construction of combinatorial bases of standard modules of affine Lie algebras was given in \cite{p1,p2}. In this construction a new and interesting class of subspaces of standard modules emerged---the so-called Feigin-Stoyanovsky's \textit{type} subspaces which coincide with principal subspaces only for affine Lie algebra $\mathfrak{sl}(2,\mathbb{C})^{\widetilde{}}$. It turned out that in the case of affine Lie algebra $\mathfrak{sl}(\ell + 1, \mathbb{C})^{\widetilde{}}$ these combinatorial bases are parametrized by $(k, \ell + 1)$-admissible configurations, combinatorial objects introduced and studied in \cite{FJLMM} and \cite{FJMMT}. In the case of other ``classical'' affine Lie algebras only bases of basic modules are constructed by using the crystal base \cite{KKMMNN} character formula.

Inspired by the use of intertwining operators in the work of Capparelli, Lepowsky and Milas, in \cite{p3} a simpler proof for the existence of combinatorial bases of Feigin-Stoyanovsky's type subspaces was given in the case of affine Lie algebra $\mathfrak{sl}(2,\mathbb{C})^{\widetilde{}}$. In this paper we extend this approach to obtain exact sequences of Feigin-Stoyanovsky's type subspaces at fixed level $k$ for affine Lie algebra $\mathfrak{sl}(l+1,\mathbb{C})^{\widetilde{}}$. In order to state
the main result we need some notation.

Denote by $\mathfrak{g} = \mathfrak{sl}(\ell + 1,\mathbb{C})$ and let $\mathfrak{h}$ be a Cartan subalgebra of $\mathfrak{g}$ with the corresponding root system $R$. Let \begin{align*} \mathfrak{g} = \mathfrak{h} + \sum_{\alpha \in R} \mathfrak{g}_{\alpha}
\end{align*} be the root space decomposition of $\mathfrak{g}$ with fixed root vectors $x_{\alpha}$. Let \begin{align*} \mathfrak{g} = \mathfrak{g}_{-1} \oplus \mathfrak{g}_0 \oplus \mathfrak{g}_1\end{align*} be a chosen $\mathbb{Z}$-grading such that $\mathfrak{h} \subset \mathfrak{g}_0$. Grading on $\mathfrak{g}$ induces $\mathbb{Z}$-grading on $\tilde{\mathfrak{g}}$: \begin{align*} \tilde{\mathfrak{g}} = \tilde{\mathfrak{g}}_{-1} \oplus \tilde{\mathfrak{g}}_{0} \oplus \tilde{\mathfrak{g}}_{1}, \end{align*} where $\tilde{\mathfrak{g}}_{1} = \mathfrak{g}_1 \otimes \mathbb{C}[t,t^{-1}]$ is a commutative Lie algebra with a basis given by \begin{align*} \{ x_{\gamma}(j) \mid j \in \mathbb{Z}, \gamma \in \Gamma \}. \end{align*} Here $\Gamma$ denotes the set of roots belonging to $\mathfrak{g}_1$.

For a standard $\tilde{\mathfrak{g}}$-module of level $k=\Lambda(c)$ with a highest weight vector $v_{\Lambda}$ define Feigin-Stoyanovsky's type subspace $W(\Lambda)$ as the following subspace of $L(\Lambda)$: \begin{align*} W(\Lambda) = U(\tilde{\mathfrak{g}}_1)\cdot v_{\Lambda}. \end{align*}

By using a description of combinatorial bases for Feigin-Stoyanovsky's type subspaces for affine Lie algebra $\mathfrak{sl}(l+1,\mathbb{C})^{\widetilde{}}$ in terms of $(k, \ell + 1)$-admissible vectors, as well as operators $\varphi_0, \dots, \varphi_{m-1}$ (constructed with certain intertwining operators between standard modules) and a simple current operator $[\omega]$, we obtain exact sequences of these subspaces at fixed level $k$ (cf. Theorem \ref{exactness_1}): \begin{align*} 0&\rightarrow W_{k_{\ell},k_0,k_1, \dots, k_{\ell-1}} \xrightarrow{[\omega]^{\otimes k}} W \xrightarrow{\varphi_0} \sum_{I_1 \in D_1(K)} W_{I_1} \xrightarrow{\varphi_1} \dots \xrightarrow{\varphi_{m-1}} W_{I_m} \rightarrow 0. \end{align*} Here $W = W(k_0 \Lambda_0 + \dots + k_{\ell} \Lambda_{\ell})$ is a Feigin-Stoyanovsky's type subspace for level $k = k_0 + \dots + k_{\ell}$, $W_{k_{\ell},k_0,k_1, \dots, k_{\ell-1}} = W(k_{\ell} \Lambda_0 + k_0 \Lambda_1 + \dots + k_{\ell - 1} \Lambda_{\ell})$, and $W_{I_1}, \dots, W_{I_m}$ Feigin-Stoyanovsky's type subspaces derived from $W$ by a procedure described in Subsection 5.1.

 Exact sequences described above directly lead to systems of relations among formal characters of those subspaces (cf. equation \ref{recurrence_short}): \begin{align*} & \sum_{I \in D(K)} (-1)^{|I|} \chi(W_I)(z_1,\dots, z_{\ell};q) = \\ & = (z_1q)^{k_0} \dots (z_{\ell}q)^{k_{\ell-1}} \chi(W_{k_{\ell},k_0,\dots, k_{\ell-1}})(z_1q,\dots, z_{\ell}q;q). \end{align*}

The paper is organized as follows. Section 2 gives the setting. In Section 3 we define Feigin-Stoyanovsky's type subspaces and present the result on combinatorial bases. Section 4 gives the vertex operator construction for fundamental modules, and introduces the intertwining operators and a simple current operator. The last two sections contain the main results of the paper: Section 5 states the result on exactness (cf. Theorem \ref{exactness_1}), while in Section 6 we obtain the corresponding system of relations among characters of Feigin-Stoyanovsky's type subspaces at fixed level (cf. equations \ref{recurrence_short} and \ref{A_recurrence}), and present the proof that such system has a unique solution.

I sincerely thank Mirko Primc for his valuable suggestions.

\section{Affine Lie algebra $\mathfrak{sl}(l+1,\mathbb{C})^{\widetilde{}}$ and standard modules}

Let $\mathfrak{g}=\mathfrak{sl}(\ell+1,\mathbb{C})$, $\ell \in \mathbb{N}$, and $\mathfrak{h}$ a Cartan subalgebra of $\mathfrak{g}$. Denote by $R$ the corresponding root system (identified in the usual way as a subset of $\mathbb{R}^{\ell +1}$): \begin{align*}R=\{\pm (\epsilon_i - \epsilon_j) \mid 1 \leq i < j \leq \ell+1 \}.\end{align*} As usual, we fix simple roots $\alpha_i=\epsilon_i - \epsilon_{i+1}$, $i=1, \dots , \ell$, and have triangular decomposition $\mathfrak{g} = \mathfrak{n}_{-} \oplus \mathfrak{h} \oplus \mathfrak{n}_{+}$. Denote by $Q = Q(R)$ the root lattice and by $P = P(R)$ the weight lattice of $R$. Let $\omega_i$, $i=1,\dots,\ell$, denote the corresponding fundamental weights. Define also $\omega_0:=0$ for later purposes. Let $\langle \cdot,\cdot\rangle$ be the Killing form on $\mathfrak{g}$. Identify $\mathfrak{h}$ and $\mathfrak{h}^{*}$ using this form, denote by $x_{\alpha}$ fixed root vectors, and normalize the form in such a way that for the maximal root $\theta$ holds $\langle\theta, \theta\rangle = 2$.

Associate to $\mathfrak{g}$ the affine Lie algebra $\tilde{\mathfrak{g}}$ (cf. \cite{kac}) \begin{align*}\tilde{\mathfrak{g}} = \mathfrak{g} \otimes \mathbb{C}[t,t^{-1}] \oplus \mathbb{C} c \oplus \mathbb{C} d,\end{align*} with Lie product given by \begin{align*}[x\otimes t^m, y \otimes t^n] = [x,y]\otimes t^{m+n} + m \langle x,y\rangle\delta_{m+n,0}c, \end{align*} for $x,y \in \mathfrak{g}$, $m,n \in \mathbb{Z}$. Here $c$ denotes the canonical central element, and $d$ the degree operator: $[d,x \otimes t^n] = nx\otimes t^n$.

We denote $x(n)=x\otimes t^n$ for $x \in \mathfrak{g}$, $n \in \mathbb{Z}$, and define the following formal Laurent series in formal variable $z$: \begin{align*}x(z)=\sum_{n\in \mathbb{Z}} x(n)z^{-n-1}.\end{align*}

Furthermore, having denoted $\mathfrak{h}^{e} = \mathfrak{h} \oplus \mathbb{C}c \oplus \mathbb{C}d$, $\tilde{\mathfrak{n}}_{\pm} = \mathfrak{g}\otimes t^{\pm 1} \mathbb{C}[t^{\pm 1}]\oplus \mathfrak{n}_{\pm}$, we have the triangular decomposition for $\tilde{\mathfrak{g}}$: \begin{align*} \tilde{\mathfrak{g}} = \tilde{\mathfrak{n}}_{-} \oplus \mathfrak{h}^{e} \oplus \tilde{\mathfrak{n}}_{+}. \end{align*} As usual, denote by $\{ \alpha_0, \alpha_1, \dots, \alpha_{\ell}\} \subset (\mathfrak{h}^{e})^{*}$ the corresponding set of simple roots, and by $\Lambda_0, \Lambda_1, \dots, \Lambda_{\ell}$ fundamental weights.

Let $L(\Lambda)$ denote standard (i.e., integrable highest weight) $\tilde{\mathfrak{g}}$-module with dominant integral highest weight \begin{align*}\Lambda = k_0\Lambda_0 + k_1 \Lambda_1 + \cdots + k_{\ell} \Lambda_{\ell},\quad k_0, k_1, \dots, k_{\ell} \in \mathbb{Z}_{+}.\end{align*} Define the level of $L(\Lambda)$ as $k=\Lambda(c)=k_0+k_1+\cdots + k_{\ell}$.

\section{Feigin-Stoyanovsky's type subspaces}

\subsection{Definition}

Fix the minuscule weight $\omega = \omega_{\ell}$ and define the following alternative basis for $\mathfrak{h}^{*}$: \begin{align*}\Gamma= \{ \alpha \in R \mid \omega(\alpha)=1 \} = \{ \gamma_1, \gamma_2, \dots , \gamma_{\ell} \mid \gamma_i=\epsilon_i - \epsilon_{\ell+1} = \alpha_i + \dots + \alpha_{\ell} \}.\end{align*} Consequently we obtain $\mathbb{Z}$-grading of $\mathfrak{g}$: \begin{align*}\mathfrak{g} = \mathfrak{g}_{-1} + \mathfrak{g}_0 + \mathfrak{g}_1,\end{align*} where $\mathfrak{g}_0 = \mathfrak{h} + \sum_{\omega(\alpha)=0} \mathfrak{g}_{\alpha}$, $\mathfrak{g}_{\pm 1} = \sum_{\alpha \in \pm \Gamma} \mathfrak{g}_{\alpha}$. The corresponding $\mathbb{Z}$-grading on affine Lie algebra $\tilde{\mathfrak{g}}$ is \begin{align*} \tilde{\mathfrak{g}} = \tilde{\mathfrak{g}}_{-1} + \tilde{\mathfrak{g}}_0 + \tilde{\mathfrak{g}}_1,\end{align*} with $\tilde{\mathfrak{g}}_0 = \mathfrak{g}_0 \otimes \mathbb{C}[t,t^{-1}]\oplus \mathbb{C} c \oplus \mathbb{C}d$, $ \tilde{\mathfrak{g}}_{\pm 1} = \mathfrak{g}_{\pm 1}\otimes\mathbb{C}[t,t^{-1}]$. Note that \begin{align*}\tilde{\mathfrak{g}}_1= \textrm{span} \{ x_{\gamma} (n) \mid \gamma \in \Gamma, n \in \mathbb{Z} \}\end{align*} is a commutative subalgebra and a $\tilde{\mathfrak{g}}_0$-module.

For an integral dominant weight $\Lambda$ define the Feigin-Stoyanovsky's type subspace of $L(\Lambda)$ as \begin{align*}W(\Lambda) = U(\tilde{\mathfrak{g}}_1) \cdot v_{\Lambda},\end{align*} $U(\tilde{\mathfrak{g}}_1)$ denoting the universal enveloping algebra for $\tilde{\mathfrak{g}}_1$, and $v_{\Lambda}$ a fixed highest weight vector of $L(\Lambda)$.

We proceed by defining colored partitions as maps \begin{align*}\pi:\{ x_{\gamma}(-j) \mid \gamma \in \Gamma, j \geq 1 \} \to \mathbb{Z}_{+}\end{align*} with finite support. For given $\pi$ define the monomial $x(\pi)$ in $U(\tilde{\mathfrak{g}}_1)$ as \begin{align*}x(\pi) = \prod x_{\gamma}(-j)^{\pi(x_{\gamma}(-j))}.\end{align*} Since every $\pi$ can be identified with a sequence $(a_i)_{i=0}^{\infty}$ with finitely many nonzero elements via $a_{\ell(j-1)+r-1}=\pi(x_{\gamma_r}(-j))$, for corresponding $x(\pi)$ we will write \begin{align}\label{monomial} x(\pi) = \dots x_{\gamma_1}(-2)^{a_{\ell}}x_{\gamma_{\ell}}(-1)^{a_{\ell-1}} \cdots x_{\gamma_1}(-1)^{a_0}.\end{align}

From Poincaré-Birkhoff-Witt theorem it follows that monomial vectors $x(\pi) v_{\Lambda}$, with $x(\pi)$ as in \eqref{monomial}, span $W(\Lambda)$.

\subsection{Combinatorial basis for $W(\Lambda)$}

For given $\Lambda = k_0 \Lambda_0 + k_1 \Lambda_1 + \dots + k_{\ell} \Lambda_{\ell}$ a monomial $x(\pi)$ given by \eqref{monomial}, or a monomial vector $x(\pi)v_{\Lambda}$ for such $x(\pi)$, is called $(k,\ell +1)$-admissible for $\Lambda$ (or $W(\Lambda)$) if the following inequalities are met: \begin{align} \nonumber a_0 &\leq k_0 \\ \label{initial} a_0+a_1 &\leq k_0 + k_1 \\ \nonumber \dots \\ \nonumber a_0 + a_1 + \dots + a_{\ell-1} &\leq k_0 + \dots + k_{\ell-1}, \end{align} and \begin{align} \label{difference} a_i + \dots a_{i+\ell} \leq k, \quad i \in \mathbb{Z}_{+}.\end{align}

We say that \eqref{initial} are the initial conditions for given $\Lambda$ (or $W(\Lambda)$), and that \eqref{difference} are the difference conditions (note that the difference conditions \eqref{difference} do not depend on the choice of $\Lambda$).

In \cite{p3} Primc demonstrated, by using certain coefficients of intertwining operators between fundamental modules as well as a simple current operator, that the spanning set for $W(\Lambda)$ given in 3.1 can be reduced to basis consisting of $(k,\ell + 1)$-admissible vectors: \begin{theorem} \label{bases} The set of $(k, \ell+1)$-admissible monomial vectors $x(\pi)v_{\Lambda}$ is a basis of $W(\Lambda)$. \end{theorem}

Note that the formulation of Theorem \ref{bases} in terms of $(k, \ell + 1)$-admissible monomials was first given in \cite{FJLMM}. We will not give details of the proof here, although some of the content needed for it is to be presented in the the next section.

\section{Intertwining operators and simple current operator}

\subsection{Fundamental modules}

In this section we state some facts needed in the following sections. For precise definitions and proofs see \cite{dl,fhl1,flm,GL,li,LL}.

Let us state the basics of the well-known vertex operator construction for fundamental $\tilde{\mathfrak{g}}$-modules (cf. \cite{fk,s}). Define Fock space $M(1)$ as induced $\hat{\mathfrak{h}}$-module \begin{align*} M(1)=U(\hat{\mathfrak{h}}) \otimes_{U(\mathfrak{h} \otimes \mathbb{C}[t] \oplus \mathbb{C}c)} \mathbb{C},\end{align*} where $\hat{\mathfrak{h}} = \mathfrak{h}\otimes \mathbb{C}[t,t^{-1}] \oplus \mathbb{C} $, such that $\mathfrak{h} \otimes \mathbb{C}[t]$ acts trivially and $c$ as identity on one-dimensional module $\mathbb{C}$. Denote by $\{ e^{\lambda} \mid \lambda \in P\}$ a basis of the group algebra $\mathbb{C}[P]$ of $P$. Then $M(1) \otimes \mathbb{C}[P]$ is an $\hat{\mathfrak{h}}$-module: $\hat{\mathfrak{h}}_{\mathbb{Z}} = \coprod_{n \in \mathbb{Z} \backslash \{ 0 \}} \mathfrak{h} \otimes t^n \oplus \mathbb{C}c$ acts as $\hat{\mathfrak{h}}_{\mathbb{Z}} \otimes 1$, and $\mathfrak{h}=\mathfrak{h} \otimes t^0$ as $1\otimes \mathfrak{h}$, with $h(0)$ given by $h(0) . e^{\lambda} = \langle h,\lambda\rangle e^{\lambda}$ for $h \in \mathfrak{h}, \lambda \in P$.

For $\lambda \in P$ we have vertex operators: \begin{align} \label{vertex} Y(e^{\lambda},z)=E^{-}(-\lambda, z)E^{+}(-\lambda, z) \epsilon_{\lambda}e^{\lambda}z^{\lambda},\end{align} with \begin{align*} E^{\pm}(\lambda, z) &= \textrm{exp} \Big( \sum_{n \geq 1} \lambda (\pm n) \frac{z^{\mp n}}{\pm n} \Big),\quad \lambda \in P \\ \epsilon_{\lambda}e^{\lambda}(v\otimes e^\mu) & = \epsilon(\lambda, \lambda + \mu)v\otimes e^{\lambda + \mu} \\ z^{\lambda}(v\otimes e^{\mu}) &= z^{\langle \lambda , \mu \rangle} v\otimes e^{\mu}, \quad v \in M(1), \lambda, \mu \in P, \end{align*} where $\epsilon$ is a 2-cocyle corresponding to a central extension of $P$ by certain finite cyclic group (cf. \cite{flm,dl}).

The action of $\hat{\mathfrak{h}}$ extends to the action of $\tilde{\mathfrak{g}}$ via \eqref{vertex} in the following manner: $x_{\alpha}(n)$ acts on $M(1)\otimes \mathbb{C}[P]$ as the coefficient of $x^{-n-1}$ in $Y(e^{\alpha},z)$: \begin{align*} x_{\alpha}(z)=Y(e^{\alpha},z),\end{align*} and $d$ as degree operator. Furthermore, the direct summands $M(1)\otimes e^{\omega_i}\mathbb{C}[Q]$ of $M(1)\otimes \mathbb{C}[P]$ are exactly standard $\tilde{\mathfrak{g}}$-modules $L(\Lambda_i)$, and we can identify highest weight vectors $v_{\Lambda_i} = 1\otimes e^{\omega_i}$, $i=0,1,\dots, \ell$.

For $\lambda \in P$ we also have Dong-Lepowsky's intertwining operators: \begin{align*}\mathcal{Y}(e^{\lambda},z):=E^{-}(-\lambda, z)E^{+}(-\lambda, z)e_{\lambda}z^{\lambda}e^{i\pi \lambda}c(\cdot, \lambda),\end{align*} where $e_{\lambda} = e^{\lambda}\epsilon(\lambda, \cdot)$, $c(\alpha, \beta)$ defined in \cite{dl}, equation (12.52).

For $\lambda_i := \omega_i - \omega_{i-1}$, $i =1, \dots, \ell$, we have: \begin{align*} [Y(e^{\gamma},z_1), \mathcal{Y}(e^{\lambda_i}, z_2)] = 0,\quad \gamma \in \Gamma, \end{align*} which implies that every coefficient of all $\mathcal{Y}(e^{\lambda_i}, z)$ commutes with $x_{\gamma}(n)$, $\gamma \in \Gamma$, $n \in \mathbb{Z}$. In particular, this holds for coefficients \begin{align*}[i]=\textrm{Res} z^{-1-\langle \lambda_i, \omega_{i-1} \rangle}c_i \mathcal{Y}(e^{\lambda_i}, z), \quad i=1,\dots, \ell,\end{align*} which we also call intertwining operators.

We have \begin{align*}L(\Lambda_0)\xrightarrow{[1]} L(\Lambda_1) \xrightarrow{[2]} L(\Lambda_2) \xrightarrow{[3]} \dots \xrightarrow{[\ell-1]} L(\Lambda_{\ell-1}) \xrightarrow{[\ell]} L(\Lambda_{\ell}),\end{align*} and for suitably chosen $c_i$ in the definition of intertwining operators also \begin{align*}v_{\Lambda_0} \xrightarrow{[1]} v_{\Lambda_1} \xrightarrow{[2]} v_{\Lambda_2} \xrightarrow{[3]} \dots \xrightarrow{[\ell-1]} v_{\Lambda_{\ell-1}} \xrightarrow{[\ell]} v_{\Lambda_{\ell}}.\end{align*}

Next, define linear bijection $[\omega]$ on $M(1)\otimes \mathbb{C}[P]$ by \begin{align*}[\omega]=e^{\omega_{\ell}}\epsilon(\cdot, \omega_{\ell}).\end{align*} We call $[\omega]$ a simple current operator (cf. \cite{DLM}). It can be shown that \begin{align*} L(\Lambda_0) \xrightarrow{[\omega]} L(\Lambda_{\ell}) \xrightarrow{[\omega]} L(\Lambda_{\ell-1}) \xrightarrow{[\omega]} \dots \xrightarrow{[\omega]} L(\Lambda_1) \xrightarrow{[\omega]} L(\Lambda_0),\end{align*} and from vertex operator formula \eqref{vertex} we have \begin{align}\label{omega1} [\omega] v_{\Lambda_0} = v_{\Lambda_{\ell}}, \quad [\omega] v_{\Lambda_i} = x_{\gamma_i}(-1)v_{\Lambda_{i-1}}, \quad i=1,\dots, \ell.\end{align} Also, from \eqref{vertex} it follows $x_{\alpha}(z)[\omega] = [\omega]z^{\langle \omega_{\ell}, \alpha \rangle} x_{\alpha}(z)$ for $\alpha \in R$ or, written in components: $x_{\alpha}(n)[\omega] = [\omega] x_{\alpha}(n + \langle \omega_{\ell}, \alpha \rangle)$, $\alpha \in R$, $n \in \mathbb{Z}$, which specially for $\gamma \in \Gamma$ gives \begin{align}\label{omega2} x_{\gamma}(n)[\omega] = [\omega] x_{\gamma}(n + 1).\end{align} In general, for every $x(\pi)v_{\Lambda}$ we have \begin{align*} [\omega] x(\pi) =x(\pi^{-}) [\omega], \end{align*} where $x(\pi^{-})$ denotes monomial obtained from $x(\pi)$ by lowering the degree of every constituting factor in $x(\pi)$ by one.

\subsection{Higher level standard modules}

Because of complete reducibility of tensor product of standard modules we can embed level $k$ standard $\tilde{\mathfrak{g}}$-module $L(\Lambda)$ with $\Lambda = k_0 \Lambda_0 + \dots + k_{\ell} \Lambda_{\ell}$ in the appropriate $k$-fold tensor product of fundamental modules \begin{align*}L(\Lambda)\subset L(\Lambda_{\ell})^{\otimes k_{\ell}} \otimes \cdots \otimes L(\Lambda_1)^{\otimes k_1} \otimes L(\Lambda_0)^{\otimes k_0},\end{align*} and we can take the corresponding highest weight vector $v_{\Lambda}$ to be \begin{align*}v_{\Lambda} = v_{\Lambda_{\ell}}^{\otimes k_{\ell}} \otimes \cdots \otimes v_{\Lambda_1}^{\otimes k_1} \otimes v_{\Lambda_0}^{\otimes k_0}.\end{align*}

Denote again $[i]:=1\otimes \cdots 1 \otimes [i] \otimes 1 \cdots \otimes 1$, $i=1,\dots, \ell$, i.e. the $k$-fold tensor product of $[i]$ with identity maps: \begin{align} \label{intertwining} & [i]: L(\Lambda_{\ell})^{\otimes k_{\ell}} \otimes \dots \otimes L(\Lambda_i)^{\otimes k_i} \otimes L(\Lambda_{i-1})^{\otimes k_{i-1}} \otimes \dots \otimes L(\Lambda_0)^{\otimes k_0} \to \\ \nonumber & \to L(\Lambda_{\ell})^{\otimes k_{\ell}} \otimes \dots \otimes L(\Lambda_i)^{\otimes k_i+1} \otimes L(\Lambda_{i-1})^{\otimes k_{i-1}-1} \otimes \dots \otimes L(\Lambda_0)^{\otimes k_0} \end{align} for $k_i \geq 1$. These are again linear maps between corresponding standard level $k$ $\tilde{g}$-modules that map highest weight vector into highest weight vector. Also, these maps commute with action of $x_{\gamma}(n)$, $\gamma \in \Gamma$, $n \in \mathbb{Z}$.

On $k$-fold tensor products of standard $\tilde{\mathfrak{g}}$-modules we use also $[\omega]^{\otimes k}$ (we will denote it again $[\omega]$), a linear bijection for which the commutation formula analogous to \eqref{omega2} holds.

\section{Exact sequences of Feigin-Stoyanovsky's type subspaces}

\subsection{Exactness and supplementary results}

For fixed level $k$ and arbitrary nonnegative integers $k_0,\dots,k_{\ell}$ such that $k_0 + \dots + k_\ell = k$ denote \begin{align*}W_{k_0, k_1, \dots , k_{\ell}} &= W(k_0\Lambda_0 + k_1 \Lambda_1 + \dots + k_{\ell} \Lambda_{\ell}), \\v_{k_0, k_1, \dots , k_{\ell}}&= v_{k_0\Lambda_0 + k_1 \Lambda_1 + \dots + k_{\ell} \Lambda_{\ell}}=v_{\ell}^{\otimes k_{\ell}} \otimes \dots \otimes v_1^{\otimes k_1}\otimes v_0^{\otimes k_0}.\end{align*} Also, let $\mathcal{B}_{k_0,k_1,\dots, k_{\ell}}$ be the set of $(k, \ell + 1)$-admissible monomials $x(\pi)$ for $k_0\Lambda_0 + k_1 \Lambda_1 + \dots + k_{\ell} \Lambda_{\ell}$. Then vectors $x(\pi)v_{k_0,k_1,\dots, k_{\ell}}$ belong to basis of $W_{k_0,k_1,\dots, k_{\ell}}$.

Let us now fix some $K=(k_0, \dots, k_{\ell})$ satisfying $k_0+\dots + k_{\ell} = k$ and define $m=\sharp \{ i=0,\dots, \ell-1 \mid k_i \neq 0 \}$. Denote \begin{align*}W=W_{k_0,k_1,\dots, k_{\ell}}, v= v_{k_0,k_1, \dots, k_{\ell}},\mathcal{B} = \mathcal{B}_{k_0,k_1,\dots, k_{\ell}}.\end{align*} For $t \in \{ 0, \dots, m-1 \}$ define \begin{align*} D_{t+1}(K) = \{ \{ i_0, \dots, i_t \} \mid 0 \leq i_0 < \dots < i_t \leq \ell-1 \ \textrm{such that} \ k_{i_j} \neq 0, j=0,\dots, t \}. \end{align*} Let $D_0(K):=\{ \emptyset \}$ and denote by $D(K)$ the union of all $D_0(K), D_1(K), \dots, D_m(K)$ (note that $D_m(K) = \{ I_m \}$ is monadic).

For $I_{t+1} = \{ i_0, \dots, i_t \} \in D_{t+1}(K)$, $t=0, \dots, m-1$, we introduce the following notation: \begin{align*}W_{I_{t+1}} = W_{\{i_0, \dots, i_t \}} = W_{k_0, \dots, k_{i_0}-1, k_{i_0+1}+1, \dots, k_{i_t}-1, k_{i_t + 1}+1, \dots, k_{\ell}},\end{align*} similarly also for highest weight vector $v_{\{ i_0, \dots , i_t \}}$ and $\mathcal{B}_{\{i_0, \dots, i_t \}}$. Also, define $W_{\emptyset} = W$, $v_{\emptyset} = v$, $\mathcal{B}_{\emptyset} = \mathcal{B}$.

By using $[i]$ given by \eqref{intertwining}, we construct the following mappings that commute with the action of $x_{\gamma}(n)$, $\gamma \in \Gamma$, $n \in \mathbb{Z}$: for $t=0, \dots , m-1$, let $\varphi_t $ be a $U(\tilde{\mathfrak{g}}_1)$-homogeneous mapping \begin{align*}\varphi_t : \sum_{I_t \in D_t(K)} W_{I_t} \to \sum_{I_{t+1} \in D_{t+1}(K)} W_{I_{t+1}}\end{align*} given component-wise by \begin{align*}\varphi_t(v_{I_t}) = \sum_{\genfrac{}{}{0pt}{}{\{ i \} \in D_1(K)}{i \notin I_t}} (-1)^{p_{I_t}(i)} v_{I_t \cup \{ i \}},\end{align*} with $p_{I_t}(i)$ such that $i=j_{p_{I_t}(i)}$ in $I_t \cup \{ i \} = \{ j_0, \dots, j_t\}$ for some $\{ j_0, \dots, j_t\} \in D_{t+1}(K)$.

Observe that $\varphi_t$ can be described also by providing the components of the image: for $I_{t+1} = \{ i_0, \dots, i_t \} \in D_{t+1}(K)$ we have \begin{align*} \varphi_t(w)_{I_{t+1}} = \sum_{0 \leq s \leq t} (-1)^s a(w)_{I_{t+1} \setminus \{ i_s\}} v_{I_{t+1}},\end{align*} where $a(w)_{I_t}$ stands for the monomial part of $w = \sum_{I_t \in D_t(K)} a(w)_{I_t} v_{I_t}$ in $W_{I_t}$.

We can now state the exactness result:

\begin{theorem}\label{exactness_1} For every $K=(k_0,\dots, k_{\ell})$ such that $k_0+ \dots + k_{\ell} = k$ the following sequence is exact: \begin{align*} 0&\rightarrow W_{k_{\ell},k_0,k_1, \dots, k_{\ell-1}} \xrightarrow{[\omega]^{\otimes k}} W \xrightarrow{\varphi_0} \sum_{I_1 \in D_1(K)} W_{I_1} \xrightarrow{\varphi_1} \dots \xrightarrow{\varphi_{m-1}} W_{I_m} \rightarrow 0. \end{align*}\end{theorem}

Note that in the case of $K=(0,0, \dots 0, k)$ we have $m=0$ which yields the short exact sequence \begin{align*}0\rightarrow W_{k,0,0, \dots, 0} \xrightarrow{ [\omega]^{\otimes k} } W_{0,0, \dots 0, k} \rightarrow 0. \end{align*}

\begin{example} For $\ell=2$, $k=2$ the exact sequences are: \begin{align*} & 0 \to W_{0,2,0} \to W_{2,0,0} \to W_{1,1,0} \to 0 \\ & 0 \to W_{0,1,1} \to W_{1,1,0} \to W_{0,2,0} \oplus W_{1,0,1} \to W_{0,1,1} \to 0 \\ & 0 \to W_{1,1,0} \to W_{1,0,1} \to W_{0,1,1} \to 0 \\ & 0 \to W_{0,0,2} \to W_{0,2,0} \to W_{0,1,1} \to 0 \\ & 0 \to W_{1,0,1} \to W_{0,1,1} \to W_{0,0,2} \to 0 \\ & 0 \to W_{2,0,0} \to W_{0,0,2} \to 0. \end{align*} \end{example}

The following lemmas (outlining some technical facts related to sets $\mathcal{B}_{A}$, $A \in D(K)$) will be used in the proof of Theorem \ref{exactness_1}.

\begin{lemma} \label{lemma1_1} For $A, B \in D(K)$ the following holds: \begin{align*} A \subset B \Rightarrow \mathcal{B}_B \subset \mathcal{B}_A. \end{align*} \end{lemma}

\begin{proof} If $A = \{ i_0, \dots, i_t \}$, the initial conditions for $W_A$ are: \begin{align} \nonumber a_0 & \leq k_0 \\ \nonumber a_0 + a_1 & \leq k_0 + k_1 \\ \nonumber & \dots \\ \nonumber a_0 + \dots + a_{i_0} &\leq k_0 + \dots + k_{i_0}-1 \\ \nonumber a_0 + \dots + a_{i_0 + 1} &\leq k_0 + \dots + k_{i_0 +1} \\ & \dots \label{lemma1_2} \\ \nonumber a_0 + \dots + a_{i_t} &\leq k_0 + \dots + k_{i_t}-1 \\ \nonumber a_0 + \dots + a_{i_t + 1} &\leq k_0 + \dots + k_{i_t +1} \\ \nonumber & \dots \\ \nonumber a_0 + \dots + a_{\ell-1} & \leq k_0 + \dots k_{\ell-1} \end{align} (in the case of $A = \emptyset$ these are just \eqref{initial}). It is obvious that the initial conditions for $W_B$ follow from those for $W_A$ if we strengthen every $(j+1)$-th inequality in \eqref{lemma1_2} for each $j \in B \setminus A$: \begin{align*}a_0 + \dots + a_j \leq k_0 + \dots + k_j -1.\end{align*} Hence, for every $x(\pi) \in \mathcal{B}_B$ also $x(\pi) \in \mathcal{B}_A$ holds. \end{proof}

\begin{lemma} \label{lemma2_1} Let $B_1, B_2 \in D(K)$. Then \begin{align*}\mathcal{B}_{B_1} \cap \mathcal{B}_{B_2} = \mathcal{B}_{B_1 \cup B_2}.\end{align*} \end{lemma}

\begin{proof} From Lemma \ref{lemma1_1} we have \begin{align*} \mathcal{B}_{B_1 \cup B_2} & \subseteq \mathcal{B}_{B_1} \\ \mathcal{B}_{B_1 \cup B_2} & \subseteq \mathcal{B}_{B_2}, \end{align*} so $\mathcal{B}_{B_1 \cup B_2} \subseteq \mathcal{B}_{B_1} \cap \mathcal{B}_{B_2}$ holds. On the other hand, a monomial $x(\pi) \in \mathcal{B}_{B_1} \cap \mathcal{B}_{B_2}$ satisfies the initial conditions both for $W_{B_1}$ and $W_{B_2}$, which means it always satisfies the stronger inequality out of two equally indexed ones, which implies $x(\pi) \in \mathcal{B}_{B_1 \cup B_2}$. \end{proof}

For $B \in D(K)$ let us define \begin{align*}\mathcal{B}^B = \mathcal{B}_B \setminus \bigcup_{\genfrac{}{}{0pt}{}{C \in D(K)}{C \supset B}} \mathcal{B}_C.\end{align*}

The next lemma describes the above defined set more explicitly:

\begin{lemma} \label{lemma3_1} For $B \in D(K)$ and $J_B = I_m \setminus B$, the set $\mathcal{B}^B$ consists of all $x(\pi)=\dots x_{\gamma_{\ell}}(-1)^{a_{\ell-1}} \dots x_{\gamma_2}(-1)^{a_1}x_{\gamma_1}(-1)^{a_0} \in \mathcal{B}_B$ such that \begin{align*}a_0 + \dots + a_{j} = k_0 + \dots + k_{j}, \quad j \in J_B.\end{align*} \end{lemma}

\begin{proof} Note that $I_m$ is the largest set in $D(K)$ and that for $B=I_m$ we have $J_B = \emptyset$. In this case the lemma trivially asserts \begin{align*}\mathcal{B}^{I_m} = \mathcal{B}_{I_m},\end{align*} which is an obvious consequence of the definition of $\mathcal{B}^{I_m}$. Next, for $B \subsetneqq I_m$ one directly checks that for every $C\supset B$, $C \in D(K)$, the set $\mathcal{B}_C$ consists of all $x(\pi) = \dots x_{\gamma_{\ell}}(-1)^{a_{\ell-1}} \dots x_{\gamma_2}(-1)^{a_1}x_{\gamma_1}(-1)^{a_0} \in \mathcal{B}_B$ that satisfy (besides the initial conditions imposed by $B$) also the following: \begin{align*}a_0 + \dots + a_i \leq k_0 + \dots + k_i -1, \quad i \in C \setminus B.\end{align*} Therefore, $\mathcal{B}^B$ consists of all $x(\pi) \in \mathcal{B}_B$ for which such inequalities do not hold (for every strict superset of $B$ belonging to $D(K)$). This means that $\mathcal{B}^B$ consists of all those $x(\pi) \in \mathcal{B}_B$ for which the following equalities take place: \begin{align*}a_0 + \dots + a_i = k_0 + \dots + k_i, \quad i \in C \setminus B, C \supset B, C \in D(K).\end{align*} All the indices from $C \in D(K)$, for all $C \supset B$, exactly comprise $J_B= I_m \setminus B$, since $I_m$ is the largest set in $D(K)$. The assertion follows. \end{proof}

\begin{lemma} \label{lemma4} For $A \in D(K)$ the set $\mathcal{B}_A$ can be expressed as a disjoint union \begin{align} \label{lemma4equation} \mathcal{B}_A = \bigcup_{\genfrac{}{}{0pt}{}{B \in D(K)}{B \supseteq A}} \mathcal{B}^B.\end{align} \end{lemma}

\begin{proof} One inclusion is clear: for $B \in D(K)$, $B \supseteq A$, we have $\mathcal{B}^B \subseteq \mathcal{B}_B \subseteq \mathcal{B}_A$. On the other hand, take $x(\pi)=\dots x_{\gamma_{\ell}}(-1)^{a_{\ell-1}} \dots x_{\gamma_2}(-1)^{a_1}x_{\gamma_1}(-1)^{a_0}$ in $\mathcal{B}_A$ and define $J = \{ j \mid a_0 + \dots + a_j = k_0 + \dots + k_j \}$. By Lemma \ref{lemma3_1} it is obvious that $x(\pi)\in \mathcal{B}^B$ for $B = I_m \setminus J$ and that $B \supseteq A$ (of course, in the case of $J=\emptyset$ we have $B= I_m$) and the equality stated in the lemma holds. Furthermore, suppose that for some $B_1, B_2 \in D(K)$, $B_1, B_2 \supseteq A$, there exists $x(\pi) \in \mathcal{B}^{B_1} \cap \mathcal{B}^{B_2}$. But, then we have $x(\pi) \in \mathcal{B}_{B_1}$ and $x(\pi) \in \mathcal{B}_{B_2}$, so $x(\pi) \in \mathcal{B}_{B_1} \cap \mathcal{B}_{B_2} = \mathcal{B}_{B_1 \cup B_2}$, which is, due to $B_1 \cup B_2 \supset B_1, B_2$, in opposition with $x(\pi) \in \mathcal{B}_{B_1}$ and $x(\pi) \in \mathcal{B}_{B_2}$. We conclude that the union on the right-hand side of \eqref{lemma4equation} must be disjoint. \end{proof}

\subsection{Proof of exactness}

We prove the Theorem \ref{exactness_1}. It is clear from the definition that $[\omega] = [\omega]^{\otimes k}$ is injective. Therefore, as the first non-trivial step we prove that $Im([\omega]^{\otimes k}) = Ker(\varphi_0)$. Recall that $\mathcal{B} = \mathcal{B}_{\emptyset}$ and let $x(\pi)v \in Ker(\varphi_0)$ for $x(\pi)=\dots x_{\gamma_{\ell}}(-1)^{a_{\ell-1}} \dots x_{\gamma_2}(-1)^{a_1}x_{\gamma_1}(-1)^{a_0} \in \mathcal{B}_{\emptyset}$. From \begin{align*}\varphi_0(x(\pi)v)= \sum_{I_1 \in D_1(K)} x(\pi)v_{I_1} =0,\end{align*} by using the fact that for $x(\pi)$ which satisfies the difference conditions and does not satisfy the initial conditions for $W(\Lambda)$ the assertion $x(\pi)v_{\Lambda} = 0$ holds, one gets \begin{align*}x(\pi) \in \mathcal{B}_{\emptyset} \setminus \bigcup_{I_1 \in D_1(K)} \mathcal{B}_{I_1} = \mathcal{B}_{\emptyset} \setminus \bigcup_{C \in D(K)} \mathcal{B}_C = \mathcal{B}^{\emptyset}.\end{align*} By using Lemma \ref{lemma3_1} we conclude that for $x(\pi)$ \begin{align*}a_0 + \dots + a_j = k_0 + \dots + k_j, \quad j \in J_{\emptyset}=I_m\end{align*} holds. Since every $k_i, i \in \{ 0, \dots, \ell-1 \} \setminus J_{\emptyset},$ equals zero, the above equalities imply equality between \textit{all} corresponding left-hand and right-hand partial sums: \begin{align*}a_0 + \dots + a_i= k_0 + \dots + k_i, \quad i \in \{ 0, \dots, \ell-1\}.\end{align*} We now have $a_i = k_i, i \in \{ 0, \dots, \ell-1 \}$, and \begin{align*} Ker(\varphi_0) = \textrm{span} \{x(\pi)v \in W \mid x(\pi)= \dots x_{\gamma_{\ell}}(-1)^{k_{\ell-1}} \dots x_{\gamma_1}(-1)^{k_0} \}.\end{align*}

On the other hand, take $x(\pi_1)v_{k_{\ell},k_0,\dots, k_{\ell-1}}$ where \begin{align*} x(\pi_1) = \dots x_{\gamma_{\ell}}(-1)^{b_{\ell-1}} \dots x_{\gamma_2}(-1)^{b_1}x_{\gamma_1}(-1)^{b_0} \in \mathcal{B}_{k_{\ell},k_0,\dots, k_{\ell-1}} \end{align*} and calculate using \eqref{omega1} and \eqref{omega2}: \begin{align*} & [\omega]^{\otimes k}(x(\pi_1)v_{k_{\ell},k_0,\dots, k_{\ell-1}}) = \\ &= x(\pi_1^{-}) x_{\gamma_{\ell}}(-1)^{k_{\ell-1}} \dots x_{\gamma_2}(-1)^{k_1} x_{\gamma_1}(-1)^{k_0}[w]^{\otimes k}(v_{k_{\ell},k_0,\dots, k_{\ell-1}})= \\ & = x(\pi_1^{-}) x_{\gamma_{\ell}}(-1)^{k_{\ell-1}} \dots x_{\gamma_2}(-1)^{k_1} x_{\gamma_1}(-1)^{k_0}v. \end{align*} If vector $x(\pi_1)$ satisfies the initial conditions for $W_{k_{\ell}, k_0, \dots, k_{\ell - 1}}$, then \begin{align*} x(\pi_1^{-}) x_{\gamma_{\ell}}(-1)^{k_{\ell - 1}} \dots x_{\gamma_2}(-1)^{k_1} x_{\gamma_1}(-1)^{k_0} v \end{align*} satisfies the difference conditions. It is therefore enough to check: \begin{align*} & b_0 \leq k_{\ell} \Rightarrow k_0 + \dots + k_{\ell-1} + b_0 \leq k \\ & b_0 + b_1 \leq k_{\ell} + k_0 \Rightarrow k_1 + \dots + k_{\ell-1} + b_0 + b_1 \leq k \\ & \dots \\ & b_0 + \dots + b_{\ell-1} \leq k_{\ell} + k_0 + \dots + k_{\ell-2} \Rightarrow k_{\ell-1} + b_0 + \dots + b_{\ell-1} \leq k.\end{align*} Also, the image of $x(\pi_1)v_{k_{\ell},k_0,\dots, k_{\ell-1}}$ obviously satisfies the initial conditions for $W$ and thus we obtain \begin{align*}Im([\omega]^{\otimes k}) = \textrm{span} \{x(\pi)v \in W \mid x(\pi)= \dots x_{\gamma_{\ell}}(-1)^{k_{\ell-1}} \dots x_{\gamma_1}(-1)^{k_0}\}.\end{align*} The assertion that $Im([\omega]^{\otimes k}) = Ker(\varphi_0)$ now immediately follows.

We continue the proof by showing $Im(\varphi_t) = Ker(\varphi_{t+1})$, $t=0, \dots, m-2$. First we show that $Im(\varphi_t) \subseteq Ker(\varphi_{t+1})$ by proving that \begin{align*}\varphi_{t+1}(\varphi_t(w))=0\end{align*} for every \begin{align*}w = \sum_{I_t \in D_t(K)} w_{I_t} \in \sum_{I_t \in D_t(K)} W_{I_t},\end{align*} $t=0, \dots, m-2$. For all $I_{t+2} = \{ i_0, \dots, i_{t+1} \} \in D_{t+2}(K)$ the following holds: \begin{align*} & \varphi_{t+1}(\varphi_t(w))_{I_{t+2}} = \\ & = \sum_{0\leq s_1 \leq t+1} (-1)^{s_1} a(\varphi_t(w))_{I_{t+2} \setminus \{ i_{s_1} \}}v_{I_{t+2}} = \\ &= \sum_{0 \leq s_2 < s_1 \leq t+1} (-1)^{s_1}(-1)^{s_2} a(w)_{I_{t+2} \setminus \{ i_{s_1}, i_{s_2} \}}v_{I_{t+2}} + \\ & + \sum_{0 \leq s_1 < s_2 \leq t+1} (-1)^{s_1}(-1)^{s_2-1} a(w)_{I_{t+2} \setminus \{ i_{s_1}, i_{s_2} \}}v_{I_{t+2}} = \\ &= \sum_{0 \leq s_2 < s_1 \leq t+1} (-1)^{s_1}(-1)^{s_2} a(w)_{I_{t+2} \setminus \{ i_{s_1}, i_{s_2} \}}v_{I_{t+2}} - \\ & - \sum_{0 \leq s_1 < s_2 \leq t+1} (-1)^{s_1}(-1)^{s_2} a(w)_{I_{t+2} \setminus \{ i_{s_1}, i_{s_2} \}}v_{I_{t+2}} = 0,\end{align*} hence $\varphi_{t+1}(\varphi_t(w)) = 0$ follows.

Next, we prove $Ker(\varphi_{t+1})\subseteq Im(\varphi_t)$. Take \begin{align*}w=\sum_{I_{t+1} \in D_{t+1}(K)} w_{I_{t+1}} \in \sum_{I_{t+1}\in D_{t+1}(K)} W_{I_{t+1}}\end{align*} such that $\varphi_{t+1}(w)=0$. This means that for every $I_{t+1} \in D_{t+1}(K)$ and $\{ i \} \in D_1(K)$ such that $i \notin I_{t+1}$ one has \begin{align*}\varphi_{t+1}(w)_{I_{t+1} \cup \{ i \}}=0.\end{align*} Fix $A = \{ i_0, \dots, i_t \} \in D_{t+1}(K)$ and arbitrary $\{ i \} \in D_1(K)$, $i \notin A$. We have $A \cup \{ i \} = \{ j_0, \dots, j_{t+1} \}$ for some $\{ j_0, \dots, j_{t+1} \} \in D_{t+2}(K)$, so $i=j_{p_A(i)}$ for some $p_A(i) \in \{ 0, \dots, t+1 \}$.

From $\varphi_{t+1}(w)_{A \cup \{ i \}}=0$ and \begin{align*} \varphi_{t+1}(w)_{A \cup \{ i \}} & = \sum_{0 \leq s \leq t+1} (-1)^s a(w)_{(A \cup \{ i \}) \setminus \{ j_s \}} v_{A \cup \{ i \}} = \\ & = (-1)^{p_A(i)} a(w)_A v_{A \cup \{ i \}} + \sum_{\genfrac{}{}{0pt}{}{0 \leq s \leq t+1}{s \neq p_A(i)}} (-1)^s a(w)_{(A \cup \{ i \}) \setminus \{ j_s \}} v_{A \cup \{ i \}} = 0 \end{align*} the following equation holds: \begin{align} \label{first} a(w)_A v_{A \cup \{ i \}} = \sum_{\genfrac{}{}{0pt}{}{0 \leq s \leq t+1}{s \neq p_A(i)}}(-1)^{s-1+ p_A(i)} a(w)_{(A \cup \{ i \}) \setminus \{ j_s \}} v_{A \cup \{ i \}}. \end{align}

By Lemma \ref{lemma4} we can write \begin{align*}a(w)_A = \sum_{\genfrac{}{}{0pt}{}{B \in D(K)}{B\supseteq A}} a(w)_A^B, \end{align*} where $a(w)_A^B$ denotes the part of $a(w)_A$ in $\mathcal{B}^B$. Analogous decomposition holds also for summands of the right-hand side of \eqref{first}. Introducing these expressions into \eqref{first} gives \begin{align} \sum_{\genfrac{}{}{0pt}{}{B \in D(K)}{B\supseteq A}} a(w)_A^B v_{A \cup \{ i \}} = \sum_{\genfrac{}{}{0pt}{}{ 0 \leq s \leq t+1}{s \neq p_A(i)}}(-1)^{s-1+p_A(i)} \sum_{\genfrac{}{}{0pt}{}{B \in D(K)}{B\supseteq (A \cup \{ i \}) \setminus \{ j_s \}}} a(w)_{(A \cup \{ i \}) \setminus \{ j_s \}}^B v_{A \cup \{ i \}}. \label{second}\end{align}

Some of the summands in \eqref{second} are trivial. We want to see for which $x(\pi) \in \mathcal{B}^B$ ($B \in D(K)$ such that $B \supseteq A$) also $x(\pi) \in \mathcal{B}_{A \cup \{ i \}}$ holds. There are two cases: \begin{itemize} \item[1.] $A \cup \{ i \} \subseteq B$: because of Lemma \ref{lemma1_1} the statement obviously holds. \item[2.] $A \cup \{ i \} \nsubseteq B$: let $x(\pi) \in \mathcal{B}^B$ and let us suppose $x(\pi) \in \mathcal{B}_{A \cup \{ i \}}$. Then $x(\pi) \in \mathcal{B}_B \cap \mathcal{B}_{A \cup \{ i \}} = \mathcal{B}_{A \cup \{ i \} \cup B} $. Since $A \cup \{ i \} \cup B \supset B$, this is in opposition to $x(\pi) \in \mathcal{B}^B$, and thus $x(\pi) \notin \mathcal{B}_{A \cup \{ i \}}$. \end{itemize}

Now \eqref{second} becomes \begin{align*} \sum_{\genfrac{}{}{0pt}{}{B \in D(K)}{B\supseteq A \cup \{ i \}}} a(w)_A^B v_{A \cup \{ i \}} = \sum_{\genfrac{}{}{0pt}{}{0 \leq s \leq t+1}{s \neq p_A(i)}} (-1)^{s-1+p_A(i)} \sum_{\genfrac{}{}{0pt}{}{B \in D(K)}{B\supseteq A \cup \{ i \}}} a(w)_{(A \cup \{ i \}) \setminus \{ j_s \}}^B v_{A \cup \{ i \}}, \end{align*} and for every $B \in D(K)$, $B \supseteq A$ \begin{align} a(w)_A^B = \sum_{\genfrac{}{}{0pt}{}{ 0 \leq s \leq t+1}{s \neq p_A(i)}} (-1)^{s-1+p_A(i)} a(w)_{(A \cup \{ i \}) \setminus \{ j_s \}}^B \label{fourth} \end{align} follows.

We want to show the existence of \begin{align*} z=\sum_{I_t \in D_t(K)} z_{I_t} \in \sum_{I_t \in D_t(K)} W_{I_t} \end{align*} such that $\varphi_t(z)=w,$ which implies that for every $I_{t+1} \in D_{t+1}(K)$ \begin{align}\varphi_t(z)_{I_{t+1}} = a(w)_{I_{t+1}}v_{I_{t+1}}\label{z_goes_to_w} \end{align} holds.

For arbitrary $I_t \in D_t(K)$ let us define \begin{align*} a(z)_{I_t} = \sum_{\genfrac{}{}{0pt}{}{\{ i \} \in D_1(K)}{i \notin I_t}} \sum_{\genfrac{}{}{0pt}{}{B \in D(K)}{B \supseteq I_t \cup \{ i \}}} \frac{(-1)^{p_{I_t}(i)}}{|B|} a(w)_{I_t \cup \{ i\}}^B, \end{align*} for $p_{I_t}(i)$ such that $i=j'_{p_{I_t}(i)}$ in $I_t \cup \{ i \} = \{ j'_0, \dots, j'_t \}$ for some $\{ j'_0, \dots, j'_t \} \in D_{t+1}(K)$.

Let us prove \eqref{z_goes_to_w} for $A \in D_{t+1}(K)$ we fixed before. We calculate: \begin{align} \label{calculation} \varphi_t(z)_A &= \sum_{0 \leq r \leq t} (-1)^r a(z)_{A \setminus \{ i_r \}} v_A = \\ \nonumber &= \sum_{0 \leq r \leq t} (-1)^r \sum_{\genfrac{}{}{0pt}{}{ \{ i \} \in D_1(K)}{i \notin A \setminus \{ i_r \}}} \sum_{\genfrac{}{}{0pt}{}{ B \in D(K)}{ B \supseteq (A \setminus \{ i_r \}) \cup \{ i \}}} \frac{(-1)^{p_{A \setminus \{ i_r \}}(i)}}{|B|} a(w)_{(A \setminus \{ i_r \}) \cup \{ i \}}^B v_A \end{align} and conclude (similarly as before) that nontrivial summands are only those indexed by $B \in D(K)$, $ B \supseteq (A \setminus \{ i_r \}) \cup \{ i \}$, for which also $B \supseteq A$ holds. We discuss two possibilities: if $i = i_r$, then $B \supseteq A$ and no additional conditions apply. But, in the case of $i \neq i_r$ we must request $i_r \in B$ in order for $B \supseteq A$ to be true (and then we have $B \supseteq A \cup \{ i \}$).

Therefore \eqref{calculation} can be rewritten as follows: \begin{align*} \varphi_t(z)_A &= \sum_{0 \leq r \leq t} (-1)^r \sum_{\genfrac{}{}{0pt}{}{ B \in D(K)}{B \supseteq A}} \frac{(-1)^r}{|B|} a(w)_{A}^B v_A + \\ &+ \sum_{0 \leq r \leq t} (-1)^r \sum_{\genfrac{}{}{0pt}{}{ \{i \} \in D_1(K)}{i \notin A}} \sum_{\genfrac{}{}{0pt}{}{ B \in D(K)}{B \supseteq A \cup \{ i \}}} \frac{(-1)^{p_{A \setminus \{ i_r \}}(i)}}{|B|} a(w)_{(A \setminus \{ i_r \}) \cup \{ i \}}^B v_A = \\ &= \sum_{\genfrac{}{}{0pt}{}{ B \in D(K)}{ B \supseteq A}} \sum_{0 \leq r \leq t} \frac{1}{|B|} a(w)_{A}^B v_A + \\ &+ \sum_{\genfrac{}{}{0pt}{}{ \{i \} \in D_1(K)}{i \notin A}} \sum_{\genfrac{}{}{0pt}{}{B \in D(K)}{ B \supseteq A \cup \{ i \}}} \frac{1}{|B|} \sum_{0 \leq r \leq t} (-1)^{r + p_{A \setminus \{ i_r \}}(i)} a(w)_{(A \setminus \{ i_r \}) \cup \{ i \}}^B v_A .\end{align*}

Let us now demonstrate that for every $\{ i \} \in D_1(K)$ such that $i \notin A$ and every $B \in D(K)$ such that $B \supseteq A \cup \{ i \}$ \begin{align} \sum_{0 \leq r \leq t} (-1)^{r + p_{A \setminus \{ i_r \}}(i)} a(w)_{(A \setminus \{ i_r \}) \cup
\{ i \}}^B = \sum_{\genfrac{}{}{0pt}{}{ 0 \leq s \leq t+1}{s \neq p_A(i)}} (-1)^{s-1+p_A(i)} a(w)_{(A \cup \{ i \}) \setminus \{ j_s \}}^B \label{final} \end{align} holds, with $\{j_0, \dots, j_{t+1} \} = A \cup \{ i \}$ for some $\{j_0, \dots, j_{t+1} \} \in D_{t+2}(K)$.

Since both sides of \eqref{final} contain $t+1$ summands, it is enough to show the equality of the corresponding summands. Depending on the choice of $0 \leq s \leq t+1 $, $s \neq p_{A}(i)$, we distinguish two cases: \begin{itemize} \item[1.] For chosen $s < p_A(i)$ the corresponding right-hand side summand equals the left-hand side summand that corresponds to $r=s$. Namely, for $r=s < p_A(i)$ we have $p_{A \setminus \{ i_r \}}(i)=p_A(i)-1$ and $i_r = j_s$, and thus \begin{align*}& r + p_{A \setminus \{ i_r \}}(i) = s + p_{A \setminus \{ i_s \}}(i)=s+ p_A(i)-1 = s-1 + p_A(i) \\ &(A \setminus \{ i_r \}) \cup \{ i \}=(A \setminus \{ i_s \}) \cup \{ i \} = (A \cup \{ i \}) \setminus \{ i_s \}= (A \cup \{ i \}) \setminus \{ j_s \}. \end{align*} \item[2.] For chosen $s > p_A(i)$ the corresponding right-hand side summand equals the left-hand side summand corresponding to $r=s-1$, because $s\geq p_A(i) + 1 \Rightarrow r \geq p_A(i)$ implies $p_{A \setminus \{ i_r \}}(i)=p_A(i)$ and $i_r = i_{s-1}=j_s$. Therefore: \begin{align*} & r + p_{A \setminus \{ i_r \}}(i) = s + p_{A \setminus \{ i_{s-1} \}}(i)=s-1+ p_A(i) \\ & (A \setminus \{ i_r \}) \cup \{ i \} = (A \setminus \{ i_{s-1} \}) \cup \{ i \} = (A \cup \{ i \}) \setminus \{ j_s \}. \end{align*} \end{itemize}

Using \eqref{fourth} and \eqref{final} we now transform $\varphi_t(z)_A$ as follows: \begin{align*} \varphi_t(z)_A &= \sum_{\genfrac{}{}{0pt}{}{ B \in D(K)}{B \supseteq A}} \frac{t+1}{|B|} a(w)_{A}^B v_A + \sum_{\genfrac{}{}{0pt}{}{\{i \} \in D_1(K)}{ i \notin A}} \sum_{\genfrac{}{}{0pt}{}{ B \in D(K)}{B \supseteq A \cup \{ i \}}} \frac{1}{|B|} a(w)_A^B v_A = \\ &= a(w)_A^A v_A + \sum_{\genfrac{}{}{0pt}{}{ B \in D(K)}{B \supset A}} \frac{t+1}{|B|} a(w)_{A}^B v_A + \sum_{\genfrac{}{}{0pt}{}{ \{i \} \in D_1(K)}{ i \notin A}} \sum_{\genfrac{}{}{0pt}{}{ B \in D(K)}{ B \supseteq A \cup \{ i \}}} \frac{1}{|B|} a(w)_A^B v_A. \end{align*} The last sum above obviously goes on all $B \in D(K)$, $B \supset A$, and every such fixed $B$ occurs exactly $|B| -|A|$ times. Thus we have \begin{align*} \varphi_t(z)_A &= a(w)_A^A v_A + \sum_{\genfrac{}{}{0pt}{}{ B \in D(K)}{B \supset A}} \frac{t+1}{|B|} a(w)_{A}^B v_A + \sum_{\genfrac{}{}{0pt}{}{ B \in D(K)}{B \supset A}} \frac{|B| - |A|}{|B|} a(w)_A^B v_A = \\ &= a(w)_A^A v_A + \sum_{\genfrac{}{}{0pt}{}{ B \in D(K)}{B \supset A}} a(w)_{A}^B v_A = \sum_{\genfrac{}{}{0pt}{}{ B \in D(K)}{B \supseteq A}} a(w)_{A}^B v_A = a(w)_A v_A, \end{align*} which finally proves $Ker(\varphi_{t+1})\subseteq Im(\varphi_t)$, $t=0, \dots, m-2$.

\section{Recurrence relations for formal characters}

\subsection{Formal characters and systems of recurrences}

Fix $k$ and $\Lambda= k_0\Lambda_0 + \dots + k_{\ell} \Lambda_{\ell}$ such that $k_0 + \dots + k_{\ell}=k$. Define the formal character $\chi(W)$ of a subspace $W = W(\Lambda)$ as \begin{align*}\chi(W)(z_1,\dots, z_{\ell}; q)= \sum \dim W^{m,n_1,\dots,n_{\ell}} q^m z_1^{n_1} \cdots z_{\ell}^{n_{\ell}},\end{align*} where $W^{m,n_1,\dots,n_{\ell}}$ denotes a subspace of $W$ spanned by monomial vectors $x(\pi)v_{\Lambda} = \dots x_{\gamma_1}(-2)^{a_{\ell}} x_{\gamma_{\ell}}(-1)^{a_{\ell-1}} \cdots x_{\gamma_1}(-1)^{a_0} v_{\Lambda} \in \mathcal{B}$ such that $x(\pi)$ is of degree $d(x(\pi)) = m$ and weight $w(x(\pi)) = n_1 \gamma_1 + \dots + n_{\ell} \gamma_{\ell}$, with $d(x(\pi))$ and $w(x(\pi))$ given by \begin{align*} d(x(\pi)) & = \sum_{k=0}^{\infty} \sum_{i = 1}^{\ell} (k+1)\cdot a_{i + k \cdot \ell -1} \\ w(x(\pi)) & = \sum_{k=0}^{\infty} \sum_{i = 1}^{\ell} \gamma_i \cdot a_{i + k \cdot \ell -1}. \end{align*}

Note that for every fixed choice of integers $n_i$, $i=1,\dots, \ell$, and $m$, such that \begin{align*} n_i &\geq k_{i-1},\quad i=1,\dots, \ell \\ m &\geq n_1+\dots + n_{\ell} + k_0+\dots + k_{\ell-1} \end{align*} and every \begin{align*} x(\pi)v_{k_{\ell},k_0,\dots, k_{\ell-1}} &\in W_{k_{\ell},k_0,\dots, k_{\ell-1}}^{m-(n_1+\dots + n_{\ell})-(k_0+\dots + k_{\ell-1}),n_1 - k_0,\dots,n_{\ell} - k_{\ell-1}}\end{align*} we have \begin{align*} & [\omega]^{\otimes k}(x(\pi)v_{k_{\ell},k_0,\dots, k_{\ell-1}}) = \\ & = x(\pi^{-})x_{\gamma_{\ell}}(-1)^{k_{\ell-1}}\cdots x_{\gamma_2}(-1)^{k_1} x_{\gamma_1}(-1)^{k_0}v \in W^{m,n_1,\dots,n_{\ell}}. \end{align*} Also, if we now apply mappings $\varphi_0, \varphi_1, \dots, \varphi_{m-1}$ consecutively to the above right-hand side vector, neither the degree nor the weight of monomial part of that vector will change. Taking this in consideration and by using Theorem \ref{exactness_1} we conclude that the following equality among the dimensions of weight subspaces holds: \begin{align*} & \dim W_{k_{\ell},k_0,\dots, k_{\ell-1}}^{m-(n_1+\dots + n_{\ell}) - (k_0+\dots + k_{\ell-1}),n_1 - k_0,\dots,n_{\ell} - k_{\ell-1}} - \dim W^{m,n_1,\dots, n_{\ell}} + \\ & + \sum_{I_1 \in D_1(K)} \dim W_{I_1}^{m,n_1,\dots, n_{\ell}} + \dots + (-1)^{m-1} \dim W_{I_m}^{m,n_1,\dots, n_{\ell}} = 0. \end{align*}

As an implication we have this relation among characters of Feigin-Stoyanovsky's type subspaces at fixed level $k$: \begin{align} \label{recurrence} \chi(W)(z_1,\dots, z_{\ell}; q) & - \sum_{I_1 \in D_1(K)} \chi(W_{I_1})(z_1,\dots, z_{\ell};q) + \\ \nonumber & + \sum_{I_2 \in D_2(K)} \chi(W_{I_2})(z_1,\dots, z_{\ell};q) + \\ \nonumber & + \dots + (-1)^m \chi(W_{I_m})(z_1,\dots, z_{\ell};q) = \\ \nonumber & = (z_1q)^{k_0} \dots (z_{\ell}q)^{k_{\ell-1}} \chi(W_{k_{\ell},k_0,\dots, k_{\ell-1}})(z_1q,\dots, z_{\ell}q;q). \end{align} Note that \eqref{recurrence} can be abbreviated to \begin{align} \label{recurrence_short} & \sum_{I \in D(K)} (-1)^{|I|} \chi(W_I)(z_1,\dots, z_{\ell};q) = \\ \nonumber & = (z_1q)^{k_0} \dots (z_{\ell}q)^{k_{\ell-1}} \chi(W_{k_{\ell},k_0,\dots, k_{\ell-1}})(z_1q,\dots, z_{\ell}q;q). \end{align} Since \eqref{recurrence} holds for $W = W(k_0 \Lambda_0 + \dots + k_{\ell} \Lambda_{\ell})$ for all possible choices of nonnegative integers $k_0,\dots, k_{\ell}$ such that $k_0+\dots + k_{\ell} = k$, we can say that \eqref{recurrence} is a system of relations. Our attempt will be directed towards proving the uniqueness of solution for this system.

\begin{example} In the case of $\ell=2$, $k=2$, the above system reads: \begin{align*} \chi(W_{2,0,0})(z_1,z_2;q) &= \chi(W_{1,1,0})(z_1,z_2;q) + (z_1q)^2\chi(W_{0,2,0})(z_1q,z_2q;q) \\ \chi(W_{1,1,0})(z_1,z_2;q) &= \chi(W_{0,2,0})(z_1,z_2;q) + \chi(W_{1,0,1})(z_1,z_2;q) - \\ & - \chi(W_{0,1,1})(z_1,z_2;q) + (z_1q)(z_2q)\chi(W_{0,1,1})(z_1q,z_2q;q) \\ \chi(W_{1,0,1})(z_1,z_2;q) &= \chi(W_{0,1,1})(z_1,z_2;q) + z_1q\chi(W_{1,1,0})(z_1q,z_2q;q) \\ \chi(W_{0,2,0})(z_1,z_2;q) &= \chi(W_{0,1,1})(z_1,z_2;q) + (z_2q)^2\chi(W_{0,0,2})(z_1q,z_2q;q) \\ \chi(W_{0,1,1})(z_1,z_2;q) &= \chi(W_{0,0,2})(z_1,z_2;q) + z_2q\chi(W_{1,0,1})(z_1q,z_2q;q) \\ \chi(W_{0,0,2})(z_1,z_2;q) &= \chi(W_{2,0,0})(z_1q,z_2q;q) \end{align*} \end{example}

Let us now write \begin{align}\chi(W_{k_0,\dots,k_{\ell}})(z_1,\dots, z_{\ell};q)= \sum_{n_1,\dots, n_{\ell} \geq 0} A_{k_0,\dots, k_{\ell}}^{n_1,\dots, n_{\ell}}(q)z_1^{n_1}\dots z_{\ell}^{n_{\ell}},\label{A}\end{align} where $A_{k_0,\dots, k_{\ell}}^{n_1,\dots, n_{\ell}}(q)$ denote formal series in one formal variable $q$. Putting \eqref{A} into \eqref{recurrence_short} gives us (with labels analogous to those we used in the previous section) the following system: \begin{align} \sum_{I \in D(K)} (-1)^{|I|} A_{I}^{n_1, \dots, n_{\ell}}(q) = q^{n} A_{k_{\ell},k_0,\dots, k_{\ell-1}}^{n_1-k_0, \dots, n_{\ell}-k_{\ell-1}}(q), \label{A_recurrence}\end{align} with $n= n_1+ \dots + n_{\ell}$.

In the next section we prove that the system \eqref{A_recurrence} has a unique solution, thus proving the solution for \eqref{recurrence_short} is also unique.

\subsection{Uniqueness of the solution}

For fixed $k$ we have the following result:

\begin{proposition} For all choices of nonnegative integers $k_0,\dots, k_{\ell}$ such that $k_0+\dots + k_{\ell}=k$ and all nonnegative integers $n_1,\dots, n_{\ell}$ such that $n_i\geq k_{i-1}$, $i=1,\dots, \ell$, the following assertion holds: \begin{align} A_{k_0,\dots, k_{\ell}}^{n_1,\dots, n_{\ell}}(q)=q^{n}  \!\!\!\!\! \sum_{\genfrac{}{}{0pt}{}{(a_0,\dots, a_{\ell - 1}) \in \mathcal{B}}{}}  \!\!\!\!\! A_{k-a,a_0,\dots, a_{\ell - 1}}^{n_1-a_0,\dots, n_{\ell}-a_{\ell - 1}}(q), \label{equality} \end{align} having denoted that $a= a_0 + \dots + a_{\ell - 1}$, $n= n_1+ \dots + n_{\ell}$, and that $\textbf{a}:= (a_0, \dots, a_{\ell - 1}) \in \mathcal{B}$ means that \eqref{initial} holds. \end{proposition}

\begin{proof} Index the equality \eqref{equality} by $(\ell + 1)$-tuple of lower indices of the series appearing on the left-hand side, i.e. by $(k_0,\dots, k_{\ell})$. We prove the proposition by induction on these tuples, in the usual lexicographic ordering among them.

Let us first check that \eqref{equality} is true if given for smallest such $(\ell+1)$-tuple, that is $(0,\dots, 0, k)$, which means we want to check if \begin{align*} A_{0,\dots,0, k}^{n_1,\dots, n_{\ell}}(q)=q^{n} \!\!\! \sum_{{\genfrac{}{}{0pt}{}{\textbf{a} \in \mathcal{B}_{0,\dots, 0,k}}{}}} \!\!\! A_{k-a,a_0,\dots, a_{\ell - 1}}^{n_1-a_0,\dots, n_{\ell}-a_{\ell - 1}}(q) \end{align*} holds. But, we see immediately that $\mathcal{B}_{0,\dots, 0, k} = \{ (0,\dots, 0) \}$, and therefore \eqref{equality} in this case reads \begin{align*} A_{0,\dots,0, k}^{n_1,\dots, n_{\ell}}(q) = q^{n} A_{k,0,\dots,0}^{n_1,\dots, n_{\ell}}(q),\end{align*} which is exactly the corresponding relation in \eqref{A_recurrence}, and the assertion follows.

Let us now fix some $(\ell+1)$-tuple $K=(k_0,\dots, k_{\ell})$. The equation in \eqref{A_recurrence} corresponding to $K$ is \begin{align*} \sum_{I \in D(K)} (-1)^{|I|} A_{I}^{n_1, \dots, n_{\ell}}(q) = q^{n} A_{k_{\ell},k_0,\dots, k_{\ell-1}}^{n_1-k_0, \dots, n_{\ell}-k_{\ell-1}}(q). \end{align*} Series appearing on the left-hand side indexed by $I \neq \emptyset$, $I \in D(K)$, yield $(\ell + 1)$-tuples strictly smaller than $K$. Therefore, by induction hypothesis, the proposition assertion holds for \eqref{equality} corresponding to those tuples. Thus we obtain \begin{align} \label{proof} & A^{n_1, \dots, n_{\ell}}(q)= \\ \nonumber & = q^{n} \sum_{\genfrac{}{}{0pt}{}{ I \in D(K)}{ I \neq \emptyset}} \sum_{\textbf{a} \in \mathcal{B}_{I}} (-1)^{|I|-1} A_{k-a,a_0,\dots, a_{\ell -1 }}^{n_1-a_0,\dots, n_{\ell}-a_{\ell - 1}}(q) + q^{n} A_{k_{\ell},k_0,\dots, k_{\ell-1}}^{n_1-k_0, \dots, n_{\ell}-k_{\ell-1}}(q).\end{align} Note that $\textbf{a} \in \mathcal{B}_I$, $I \in D(K) \setminus \{ \emptyset \}$, belong also to $\mathcal{B}_{\emptyset} = \mathcal{B}$ (an obvious consequence of Lemma \ref{lemma1_1}). Moreover, it is not hard to see that those $\ell$-tuples are \textit{all} $\textbf{a} \in \mathcal{B}$ except $\textbf{a} = (k_0,\dots, k_{\ell-1})$. But, the last remaining summand on the right-hand side of \eqref{proof} is exactly the one ``produced'' by this ``missing'' tuple. It only remains to show that each $\textbf{a} \in \mathcal{B}$ appears only once on the right-hand side of \eqref{proof}. For arbitrary $\textbf{a} \in \mathcal{B}$ there exists $I \in D(K)$ such that $\mathcal{B}^I$ (cf. Lemma \ref{lemma4}). But then, because of Lemma \ref{lemma1_1}, also $\textbf{a} \in B_J$ holds for all $J \in D(K)$ such that $J \subseteq I$. This means that the factor multiplying $A_{k-a,a_0,\dots, a_{\ell - 1}}^{n_1-a_0,\dots, n_{\ell}-a_{\ell - 1}}(q)$ in \eqref{proof} equals to $\sum_{j=1}^{|I|} (-1)^{j-1}{\binom{|I|}{j}}$. We show this factor equals one: \begin{align*} & 1- \sum_{j=1}^{|I|} (-1)^{j-1}{\binom{|I|}{j}} = \sum_{j=0}^{|I|} (-1)^{j}{\binom{|I|}{j}}= \sum_{j=0}^{|I|} (-1)^{j}1^{|I|-j}{\binom{|I|}{j}} = 0. \end{align*} From \eqref{proof} we now get \begin{align*} A_{k_0,\dots, k_{\ell}}^{n_1,\dots, n_{\ell}}(q)=q^{n} \sum_{\genfrac{}{}{0pt}{}{\textbf{a} \in \mathcal{B}}{}} A_{k-a,a_0,\dots, a_{\ell - 1}}^{n_1-a_0,\dots, n_{\ell}-a_{\ell - 1}}(q), \end{align*} which completes the inductive proof. \end{proof}

\begin{proposition} For all choices of nonnegative integers $k_0,\dots, k_{\ell}$ such that $k_0+\dots + k_{\ell}=k$ and all nonnegative integers $n_1,\dots, n_{\ell}$ such that $n_i\geq k_{i-1}$, $i=1,\dots, \ell$, the following assertion is true: \begin{align*} A_{k_0,\dots, k_{\ell}}^{n_1,\dots, n_{\ell}}(q) = \frac{q^n}{1-q^n} \Big( \!\!\! \sum_{\genfrac{}{}{0pt}{}{ \textbf{a} \in \mathcal{B}}{ \textbf{a} \neq (0,\dots, 0)}} \!\!\!\!\! A_{k-a,a_0,\dots, a_{\ell - 1}}^{n_1-a_0,\dots, n_{\ell}-a_{\ell - 1}}(q) + q^{n} \!\!\!\!\! \sum_{\genfrac{}{}{0pt}{}{\textbf{a} \in \mathcal{B}_{k,0,\dots, 0} \setminus \mathcal{B}}{\textbf{a} \neq (0,\dots, 0)}} \!\!\!\!\!\!\! A_{k-a,a_0,\dots, a_{\ell - 1}}^{n_1-a_0,\dots, n_{\ell}-a_{\ell - 1}}(q) \Big).\end{align*} \end{proposition}

\begin{proof} From \eqref{equality} we have \begin{align*} A_{k,0,\dots, 0}^{n_1,\dots, n_{\ell}}(q) & = q^{n} \!\!\! \sum_{\genfrac{}{}{0pt}{}{\textbf{a} \in \mathcal{B}_{k,0,\dots, 0}}{}} \!\!\! A_{k-a,a_0,\dots, a_{\ell - 1}}^{n_1-a_0,\dots, n_{\ell}-a_{\ell - 1}}(q) = \\ &= q^{n} \!\!\! \sum_{\genfrac{}{}{0pt}{}{\textbf{a} \in \mathcal{B}_{k,0,\dots, 0} }{\textbf{a} \neq (0,\dots, 0)}} \!\!\! A_{k-a,a_0,\dots, a_{\ell - 1}}^{n_1-a_0,\dots, n_{\ell}-a_{\ell - 1}}(q) + q^{n} A_{k,0,\dots, 0}^{n_1,\dots, n_{\ell}}(q), \end{align*} and it follows that \begin{align} A_{k,0,\dots, 0}^{n_1,\dots, n_{\ell}}(q) = \frac{q^n}{1-q^n} \!\! \sum_{\genfrac{}{}{0pt}{}{\textbf{a} \in \mathcal{B}_{k,0,\dots, 0} }{\textbf{a} \neq (0,\dots, 0)}} \!\!\!  A_{k-a,a_0,\dots, a_{\ell - 1}}^{n_1-a_0,\dots, n_{\ell}-a_{\ell - 1}}(q). \label{half_done} \end{align} Let us now fix nonnegative integers $k_0,\dots, k_{\ell}$ such that $k_0+\dots + k_{\ell}=k$. From \eqref{equality} by using \eqref{half_done} we get \begin{align*} & A_{k_0,\dots, k_{\ell}}^{n_1,\dots, n_{\ell}}(q) = q^{n} \sum_{\genfrac{}{}{0pt}{}{\textbf{a} \in \mathcal{B}}{}} A_{k-a,a_0,\dots, a_{\ell - 1}}^{n_1-a_0,\dots, n_{\ell}-a_{\ell - 1}}(q) = \\& = q^{n} \Big( A_{k,0,\dots, 0}^{n_1,\dots, n_{\ell}}(q) + \!\!\! \sum_{\genfrac{}{}{0pt}{}{\textbf{a} \in \mathcal{B}}{\textbf{a} \neq (0,\dots, 0)}} \!\!\! A_{k-a,a_0,\dots, a_{\ell - 1}}^{n_1-a_0,\dots, n_{\ell}-a_{\ell - 1}}(q) \Big) = \\ & = q^{n} \Big( \frac{q^{n}}{1-q^{n}} \sum_{\genfrac{}{}{0pt}{}{\textbf{a} \in \mathcal{B}_{k,0,\dots, 0} }{\textbf{a} \neq (0,\dots, 0)}} \!\!\! A_{k-a,a_0,\dots, a_{\ell - 1}}^{n_1-a_0,\dots, n_{\ell}-a_{\ell - 1}}(q) + \!\!\! \sum_{\genfrac{}{}{0pt}{}{\textbf{a} \in \mathcal{B} }{\textbf{a} \neq (0,\dots, 0)}} \!\!\! A_{k-a,a_0,\dots, a_{\ell -1}}^{n_1-a_0,\dots, n_{\ell}-a_{\ell - 1}}(q) \Big) = \\ &= \frac{q^{n}}{1-q^{n}} \Big( q^{n} \!\! \sum_{\genfrac{}{}{0pt}{}{\textbf{a} \in \mathcal{B}_{k,0,\dots, 0} }{\textbf{a} \neq (0,\dots, 0)}} \!\!\! A_{k-a,a_0,\dots, a_{\ell - 1}}^{n_1-a_0,\dots, n_{\ell}-a_{\ell - 1}}(q) + (1-q^n) \!\!\!\! \sum_{\genfrac{}{}{0pt}{}{\textbf{a} \in \mathcal{B} }{\textbf{a} \neq (0,\dots, 0)}} \!\!\! A_{k-a,a_0,\dots, a_{\ell - 1}}^{n_1-a_0,\dots, n_{\ell}-a_{\ell - 1}}(q) \Big) = \\ & = \frac{q^n}{1-q^n} \Big( \!\!\! \sum_{\genfrac{}{}{0pt}{}{\textbf{a} \in \mathcal{B}}{\textbf{a} \neq (0,\dots, 0)}} \!\!\! A_{k-a,a_0,\dots, a_{\ell - 1}}^{n_1-a_0,\dots, n_{\ell}-a_{\ell - 1}}(q) + q^{n} \!\!\! \sum_{\genfrac{}{}{0pt}{}{\textbf{a} \in \mathcal{B}_{k,0,\dots, 0} \setminus \mathcal{B}}{\textbf{a} \neq (0,\dots, 0)}} \!\!\! A_{k-a,a_0,\dots, a_{\ell - 1}}^{n_1-a_0,\dots, n_{\ell}-a_{\ell - 1}}(q) \Big), \end{align*} which proves the proposition. \end{proof}

A direct consequence of this proposition is the fact that every $A_{k_0,\dots, k_{\ell}}^{n_1,\dots, n_{\ell}}(q)$ can be expressed as a linear combination (with coefficients being rational functions of formal variable $q$) of rational functions $A_{k-a,a_0,\dots, a_{\ell - 1}}^{n_1-a_0,\dots, n_{\ell}-a_{\ell - 1}}(q)$ indexed by $\ell$-tuples $(n_1-a_0,\dots, n_{\ell}-a_{\ell - 1})$ of upper indices, which are all strictly smaller than $(n_1,\dots,n_{\ell})$, in the usual lexicographic ordering. This inductively guarantees the uniqueness of the solution for system \eqref{A_recurrence}.


\begin{thebibliography}{0000000000}
\bibitem[A1]{A1} G. Andrews, \textit{The theory of partitions}, Encyclopedia of Mathematics and its Applications, Vol. 2, Addison-Wesley, 1976.
\bibitem[A2]{A2} G. Andrews, \textit{An analytic proof of the Rogers-Ramanujan-Gordon identities}, Amer. J. Math. {\bf 88} (1966), 844--846.
\bibitem[CaLM1]{CaLM1} C. Calinescu, J. Lepowsky, A. Milas, \textit{Vertex-algebraic structure of the principal subspaces of certain $A_1^{(1)}$-modules, I: level one case}, Int. J. Math. {\bf 19} (2008), 71--92.
\bibitem[CaLM2]{CaLM2} C. Calinescu, J. Lepowsky, A. Milas, \textit{Vertex-algebraic structure of the principal subspaces of certain $A_1^{(1)}$-modules, II: higher-level case}, J. Pure Appl. Algebra, to appear
\bibitem[Ca1]{Ca1} C. Calinescu, \textit{Intertwining vertex operators and certain representations of $\widehat{\mathfrak{sl}(n)}$}, Commun. Contemp. Math. {\bf 10} (2008), 47--79.
\bibitem[Ca2]{Ca2} C. Calinescu, \textit{Principal subspaces of higher-level standard $\widehat{\mathfrak{sl}(3)}$-modules}, J. Pure Appl. Algebra {\bf 210} (2007), 559--575.
\bibitem[CLM1]{CLM1} S. Capparelli, J. Lepowsky and A. Milas, \textit{The Rogers-Ramanujan recursion and intertwining operators}, Comm. Contemp. Math. {\bf 5} (2003), 947--966.
\bibitem[CLM2]{CLM2} S. Capparelli, J. Lepowsky and A. Milas, \textit{The Rogers-Selberg recursions, the Gordon-Andrews identities and intertwining operators}, math.QA/0310080.
\bibitem[DL]{dl} C. Dong and J. Lepowsky, \textit{Generalized Vertex Algebras and Relative Vertex Operators}, Progress in Math. {\bf 112}, Birkha\"user, Boston, 1993.
\bibitem[DLM]{DLM} C. Dong, H. Li and G. Mason, \textit{Simple currents and extensions of vertex operator algebras}, Comm. Math. Physics {\bf 180} (1996), 671--707.
\bibitem[FJLMM]{FJLMM} B. Feigin, M. Jimbo, S. Loktev, T. Miwa and E. Mukhin, \textit{Bosonic formulas for $(k,\ell)$-admissible partitions}, math.QA/0107054; \textit{Addendum to `Bosonic formulas for $(k,\ell)$-admissible partitions'}, math.QA/0112104.
\bibitem[FJMMT]{FJMMT} B. Feigin, M. Jimbo, T. Miwa, E. Mukhin and Y. Takeyama, \textit{Fermionic formulas for $(k,3)$-admissible configurations}, Publ. RIMS {\bf 40} (2004), 125--162.
\bibitem[FHL]{fhl1} I. Frenkel, Y.-Z. Huang and J. Lepowsky, \textit{On axiomatic approaches to vertex operator algebras and modules}, Memoirs Amer. Math. Soc. {\bf 104}, 1993.
\bibitem[FK]{fk} I. Frenkel, V. Kac, \textit{Basic representations of affine Lie algebras and dual resonance models}, Invent. Math. {\bf 62} (1980), 23--66.
\bibitem[FLM]{flm} I. Frenkel, J. Lepowsky and A. Meurman, \textit{Vertex Operator Algebras and the Monster}, Pure and Appl. Math. {\bf 134}, Academic Press, Boston, 1988.
\bibitem[FS]{FS} A. V. Stoyanovsky and B. L. Feigin, \textit{Functional models of the representations of current algebras, and semi-infinite Schubert cells}, (Russian) Funktsional. Anal. i Prilozhen. {\bf 28} (1994), no. 1, 68--90, 96; translation in Funct. Anal. Appl. {\bf 28} (1994), no. 1, 55-72; preprint B. Feigin and A. Stoyanovsky, \textit{Quasi-particles models for the representations of Lie algebras and geometry of flag manifold}, hep-th/9308079, RIMS 942.
\bibitem[G]{G} G. Georgiev, \textit{Combinatorial constructions of modules for infinite-dimensional Lie algebras, I. Principal subspace}, J. Pure Appl. Algebra {\bf 112} (1996), 247--286.
\bibitem[GL]{GL} Y. Gao, H.-S. Li, \textit{Generalized vertex algebras generated by parafermion-like vertex operators}, J. Algebra {\bf 240} (2001), 771--807.
\bibitem[K]{kac} V.G. Kac, \textit{Infinite-dimensional Lie algebras}, 3rd ed. Cambridge University Press, Cambridge, 1990.
\bibitem[KKMMNN]{KKMMNN} S.-J. Kang, M. Kashiwara, K. C. Misra, T. Miwa, T. Nakashima and A. Nakayashiki, \textit{Affine crystals and vertex models}, International Journal of Modern Physics A, Vol. 7, Suppl. 1A, Proceedings of the RIMS Research Project 1991, ``Infinite Analysis'', World Scientific, Singapore, 1992, 449--484.
\bibitem[Li]{li} H.-S. Li, \textit{Local systems of vertex operators, vertex superalgebras and modules}, J. Pure Appl. Alg. {\bf 109} (1996), 143--195.
hep-th/9406185.
\bibitem[LL]{LL} J. Lepowsky, H.-S. Li, \textit{Introduction to Vertex Operator Algebras and Their Representations}, Progress in Math. {\bf 227}, Birkh\"auser, Boston, 2004.
\bibitem[LP]{lp} J. Lepowsky, M. Primc, \textit{Structure of the standard modules for the affine Lie Algebra $A_1^{(1)}$}, Contemporary Math. {\bf 46} (1985), 1--84.
\bibitem[LW1]{lw1} J. Lepowsky and R. L. Wilson, \textit{The structure of standard modules, I: Universal algebras and the Rogers-Ramanujan identities}, Invent. Math. {\bf 77} (1984), 199--290.
\bibitem[LW2]{lw2} J. Lepowsky and R. L. Wilson, \textit{The structure of standard modules, II: The case $A_{1}^{(1)}$, principal gradation}, Invent. Math. {\bf 79} (1985), 417--442.
\bibitem[MP]{mp} A. Meurman and M. Primc, \textit{Annihilating fields of standard modules of $\mathfrak{sl}(2,\mathbb{C})^{\widetilde{} }$ and combinatorial identities}, Memoirs Amer. Math. Soc. {\bf 652} (1999)
\bibitem[P1]{p1} M. Primc, \textit{Vertex operator construction of standard modules for $A_n^{(1)}$}, Pacific J. Math {\bf 162} (1994), 143--187.
\bibitem[P2]{p2} M. Primc, \textit{Basic Representations for classical affine Lie algebras}, J. Algebra {\bf 228} (2000), 1--50.
\bibitem[P3]{p3} M. Primc, \textit{$(k,r)$-admissible configurations and intertwining operators}, Contemp. Math. {\bf 442} (2007), 425--434.
\bibitem[S]{s} G. Segal, \textit{Unitary representations of some infinite-dimensional groups,} Commun. Math. Phys. {\bf 80} (1981), 301--342.
\end{thebibliography}
\end{document}